\documentclass[12pt]{article}

\usepackage{graphicx} 
\usepackage{amsmath} 
\usepackage{amssymb,amsthm}
\usepackage{color}

\newcommand{\R}{{\mathbb R}}\newcommand{\N}{{\mathbb N}}

\newcommand{\C}{{\mathbb C}}

\newcommand{\cA}{{\mathcal A}}
\newcommand{\cU}{{\mathcal U}}
\newcommand{\cE}{{\mathcal E}}
\newcommand{\cT}{{\mathcal T}}

\newcommand{\X}{X}
\newcommand{\T}{T}

\newcommand{\veps}{\varepsilon}

\let\epsilon\varepsilon
\let\theta\vartheta

\newtheorem{theorem}{Theorem}[section]

\newtheorem{corollary}[theorem]{Corollary}

\newtheorem{remark}[theorem]{Remark}\newtheorem{example}[theorem]{Example}

\usepackage{filecontents}

\title{Some remarks about an effective description
of high-frequency wave-packet propagation}
\author{Anna Logioti$^1$, Xin Meng$^{2}$, Guido Schneider$^1$ \\
{\small
$^1$Institut f\"ur Analysis, Dynamik und Modellierung, Universit\"at Stuttgart, } \\ {\small Pfaffenwaldring 57,  70569 Stuttgart, Germany}
\\
{\small $^{2}$School of Mathematics, Jilin University,  Changchun  130012,  People's Republic of China}
}

\begin{document}

\maketitle

\begin{abstract}
We consider systems of the form
$$
\partial_{\tau} \cU + \cA(\partial_{\xi})  \cU  + \frac{1}{\varepsilon} \cE  \cU  = \cT_{2}( \cU , \cU ) + \varepsilon \cT_3( \cU , \cU , \cU ),
$$
with $ 0 < \varepsilon \ll 1 $ a small perturbation parameter.
We are interested in an
effective description
of high-frequency wave-packet propagation
associated to highly oscillatory initial conditions
$$
 \cU (\xi,0) =  \cU_*(\xi) e^{ik_0 \xi/\varepsilon} + c.c..
$$
By classical perturbation analysis for polarized initial conditions
NLS approximations up to an arbitrary order and for non-polarized initial conditions
a system of decoupled NLS equations can be derived
for the approximate
description of the associated solutions.
Under the validity of a number of non-resonance conditions we prove error estimates
between these formal approximations and true solutions of the
original system.
The result improves results from the existing literature in at least two directions,
firstly, the handling of higher order approximations in case of quadratic nonlinearities $ \cT_2(\cU,\cU)$ and secondly,
the handling of non-polarized initial conditions.
\end{abstract}

\section{Introduction}

We consider
\begin{equation} \label{sys1}
\partial_{\tau} \cU + \cA(\partial_{\xi}) \cU + \frac{1}{\varepsilon} \cE \cU = \cT_2(\cU,\cU) + \varepsilon \cT_3(\cU,\cU,\cU),
\end{equation}
with $ \cU(\xi,\tau) \in \R^N $, $ \xi \in \R^d $, $ N,d \in \N $, $ \tau \in [0,\tau_0/\varepsilon] $ for a $ \tau_0 > 0 $,
$
\cA(\partial_{\xi})  = \sum_{j=1}^d \cA_j \partial_{\xi_j}
$,
$ \cA_j = \cA_j^T \in \R^{N\times N} $, $ \cE = -\cE^T \in \R^{N\times N} $,
$ \cT_2: \R^N \times \R^N \to \R^N $ a bilinear mapping
and $ \cT_3: \R^N \times \R^N \times \R^N \to \R^N $ a trilinear mapping.
For most of the paper we assume $ d = 1 $ due to the possible application in nonlinear optics and due to the purpose of this paper.
The initial conditions for \eqref{sys1} are given by
\begin{equation}\label{init1}
\cU(\xi,0) = \cU_*(\xi) e^{ik_0 \xi/\varepsilon} + c.c.,
\end{equation}
where $ 0 < \varepsilon \ll 1 $ is a small perturbation parameter and $ k_0 > 0 $ is a fixed spatial
wave number.
This class of problems \eqref{sys1} include, e.g., the Maxwell-Lorentz system and
Klein-Gordon systems, cf. \cite[Section 2.1]{Co02,CL09}.
This class of problems has been considered for instance in \cite{CL09,Rauch12} or recently in \cite{BJL24}
where in case $ \cT_2 = 0 $ higher order
NLS approximations have been derived and justified by establishing approximation results.
The results presented in the following improve these papers in at least two directions,
firstly, the handling of higher order NLS approximations for quadratic nonlinearities $ \cT_2(\cU,\cU)$ and secondly
the handling of more than one NLS-scaled wave-packet
which is called the non-polarized situation in the following.
We obtain these improved results by relating these questions to questions
which have already been solved in the existing literature or which can easily be generalized to the present situation,
namely first the justification of higher order NLS approximations for systems with quadratic nonlinearities and secondly the separation of internal and interaction dynamics of NLS-scaled wave-packets.

As will be clear in the following, higher order nonlinear terms
$$
\ldots + \varepsilon^2 \cT_4(\cU,\cU,\cU,\cU) + \varepsilon^3 \cT_5(\cU,\cU,\cU,\cU,\cU) + \ldots
$$
can easily be incorporated in the subsequent analysis. However, for
notational simplicity we restrict ourselves to \eqref{sys1}.

Error estimates for the NLS approximation of dispersive systems
without small perturbation parameter in the underlying equations have already been
known for a few decades.
The NLS equation as envelope equation has been derived first in \cite{Za68}.
For systems without quadratic terms a simple application of Gronwall's
inequality is sufficient to obtain such estimates \cite{KSM92}.
A very general NLS approximation result including quadratic nonlinearities
has been shown in \cite{Kal87}.
This result was improved in a number of papers by weakening the
non-resonance conditions which are necessary to eliminate the
quadratic terms in order to apply Gronwall's inequality again,
e.g. \cite{Schn98,Schn05,DHSZ16}.
An application of the theory to the water wave problem
can be found in \cite{TW12,DSW16,D21}, and quasilinear wave equations
are considered in \cite{CW18,D17,DH18}.
Another example of an NLS approximation result for dispersive systems
with a small perturbation parameter in the underlying equations exists for instance  for
the Klein-Gordon-Zakharov system from plasma physics, e.g. \cite{MN02}.
Our approach to handle higher order NLS approximations for quadratic nonlinearities $ \cT_2(\cU,\cU)$
follows \cite{Schn98}.

Our approach to handle more than one NLS-scaled wave-packet is based on the literature
about
the separation of internal and interaction dynamics of NLS-scaled wave-packets, in particular on \cite{CCSU08,CS12}.
For dispersive systems it has first been observed  in \cite{PW95}
that for spatially localized NLS-scaled wave-packets
no interaction appears in lowest order w.r.t. the small
perturbation parameter $ 0 < \varepsilon \ll 1 $.
This result has been improved in a number of papers, e.g. \cite{CSU07,CCSU08,CS12,SC15}, now allowing to separate
the internal from the interaction dynamics of the wave-packets up to high order.

The plan of the paper is as follows.
In Section \ref{secex} we present two examples  which show how
the existing theory for NLS approximations  can be made applicable
for systems of the form \eqref{sys1}.
In Section \ref{sec2} we rescale
\eqref{sys1} in such a way that the existing theory for
dispersive systems
without small perturbation parameter in the equations
becomes applicable.
In Section \ref{sec3} we explain the handling of quadratic nonlinearities $ \cT_2(\cU,\cU)$ and
in Section \ref{sec4} we explain our main result, namely 
the handling of more than one NLS-scaled wave-packet.
Finally, in Section \ref{sec6} we explain how the results change if
$ x \in \R^d $ with $ d \geq 2 $ is considered.

\medskip

{\bf Acknowledgement.} The work is partially supported by the Deutsche Forschungsgemeinschaft DFG through the SFB 1173 ''Wave phenomena'' with the Project-ID 258734477 and the China
Scholarship Council through the Project-ID  202306170135. Xin Meng is grateful to the Institute for Analysis, Dynamics,
and Modeling at the University of Stuttgart for its kind hospitality during the visit. Guido Schneider would like to thank the
Mathematische Forschungsinstitut Oberwolfach for hosting the workshop
Nonlinear Optics: Physics, Analysis, and Numerics
where this research has been started.
Moreover, he would like to thank Julian Baumstark, Tobias Jahnke and
Christian Lubich for interesting discussions.
\section{Two examples}

\label{secex}

The first example motivates our subsequent approach.
\begin{example} \label{example1}
{\rm
We consider the nonlinear Klein-Gordon equation
\begin{equation} \label{KG}
\partial_t^2 u = \partial_x^2 u - u + u^2 + u^3 ,
\end{equation}
with $ t \in \R $, $ x \in \R $, and $ u(x,t) \in \R $.
By introducing $ v = \partial_t u $ and  $ w =  \partial_x u $ we obtain the system
\begin{eqnarray*}
\partial_t u  & = & v ,\\
\partial_t v   & = & \partial_x w - u + u^2 + u^3 ,\\
\partial_t w & = & \partial_x v  .
\end{eqnarray*}
This system can be written in the form
\begin{equation} \label{sys2}
\partial_{t} U + \cA \partial_x U +  \cE U = \cT_2(U,U) +  \cT_3(U,U,U),
\end{equation}
with
$$
U = \left( \begin{array}{c} u \\ v \\ w \end{array} \right), \quad
 \cA  =  \left( \begin{array}{ccc}  0 & 0 & 0 \\ 0 & 0 & -1 \\ 0 & -1 & 0 \end{array} \right),
 \quad
 \cE  =  \left( \begin{array}{ccc}  0 & -1 & 0 \\ 1 & 0 & 0 \\ 0 & 0 & 0 \end{array} \right), \quad
$$
and
$$
\cT_2(U,U)  = \left( \begin{array}{c} 0 \\ u^2 \\ 0 \end{array} \right), \quad
  \cT_3(U,U,U) = \left( \begin{array}{c} 0 \\ u^3 \\ 0 \end{array} \right).
$$
By setting $ \tau = \varepsilon t $, $ \xi = \varepsilon x $, and $ U =  \varepsilon \cU $
this system transforms in our starting system \eqref{sys1}.

For \eqref{KG} a NLS equation
$$
\partial_T A  =  i \nu_1 \partial_X^2 A + i \nu_2 A |A|^2,
$$
with coefficients $ \nu_1 $, $  \nu_2 \in \R $,
can be derived with the ansatz
$$
u(x,t) = \varepsilon A(X,T) e^{i(k_0 x - \omega_0 t)}
+ c.c. + \mathcal{O}(\varepsilon^2),
$$
where $ X = \varepsilon (x-ct) $ and $ T = \varepsilon^2 t $.
Herein, $ k_0 $ and $ \omega_0 $ satisfy the linear dispersion relation
$ \omega_0^2 = k_0^2 + 1 $, and  $ c = \frac{d\omega_0}{dk_0} $
is the linear group velocity. By definition of $ v $ and $ w $ the associated ansatz
for the solution $ U $ of \eqref{sys2} is given by
$$
U(x,t) = \varepsilon A(X,T) \left( \begin{array}{c} 1 \\ - i \omega_0 \\ i k_0 \end{array} \right)
e^{i(k_0 x - \omega_0 t)}
+ c.c. + \mathcal{O}(\varepsilon^2).
$$
Hence,  the
associated ansatz
for the solution $ \cU $ of \eqref{sys1} is given by
$$
\cU(\xi,\tau) =  A(X,T)  \left( \begin{array}{c} 1 \\ - i \omega_0 \\ i k_0 \end{array} \right)
e^{i \varepsilon^{-1}(k_0 \xi - \omega_0 \tau)}
+ c.c. + \mathcal{O}(\varepsilon),
$$
where $ X = \xi - c \tau $, $ T = \varepsilon \tau $.
As a consequence we have for the initial condition of \eqref{sys1} that
$$
\cU(\xi,0) = A(\xi,0) \left( \begin{array}{c} 1 \\ - i \omega_0 \\ i k_0 \end{array} \right)e^{ik_0 \xi/\varepsilon} + c.c..
$$
 }\end{example}
In order to have an example for which all assumptions, used in the subsequent validity proofs, can be checked
easily, we also consider
\begin{example}
\label{example2}
{\rm
We consider
\begin{equation} \label{sys2b}
\partial_{t} U + \cA \partial_x U +  \cE U = \cT_2(U,U) +  \cT_3(U,U,U),
\end{equation}
with
$$
U = \left( \begin{array}{c} u \\ v  \end{array} \right), \quad
 \cA  =  \left( \begin{array}{ccc}   0 & 1 \\  1 & 0 \end{array} \right),
 \quad
 \cE  =  \left( \begin{array}{ccc}  0 & 1 \\  - 1 & 0 \end{array} \right), \quad
$$
and
$$
\cT_2(U,U)  = \mathcal{O}(u^2+v^2) , \quad
  \cT_3(U,U,U) = \mathcal{O}(|u|^3+|v|^3) .
$$
By setting $ \tau = \varepsilon t $, $ \xi = \varepsilon x $, and $ U =  \varepsilon \cU $
this system transforms in our starting system \eqref{sys1}.
The linear operator $ \mathcal{M}  = \cA \partial_{x}  +  \cE $ is given in Fourier space by
$$
 \widehat{\mathcal{M}}(k) = \left( \begin{array}{ccc}  0 & 1 + ik  \\  -1 + ik & 0   \end{array} \right) ,
$$
The matrix $
 \widehat{\mathcal{M}}(k)$ has two eigenvalues, namely
$$ \omega_1(k) = - i \sqrt{1+k^2}, \qquad \omega_2(k) =  i \sqrt{1+k^2} . $$}\end{example}

\section{Rescaling and transforming the system}

\label{sec2}

In this section we rescale and transform
\eqref{sys1} in such a way that the existing theory for
dispersive systems
without small perturbation parameter in the underlying equations
becomes applicable. Following Example \ref{example1} of  the last section
we introduce $ t,x,U $ by $ \tau = \varepsilon t $, $ \xi = \varepsilon x $, and $ U =  \varepsilon \cU $, and consider then
\begin{equation} \label{sys4}
\partial_{t} U +  \cA \partial_{x} U +  \cE U = \cT_2(U,U) +  \cT_3(U,U,U),
\end{equation}
with initial conditions
$$
 U (x,0) =  \varepsilon U_*(\varepsilon x ) e^{ik_0 x} + c.c.
$$
This is a special form of a dispersive wave system.
For dispersive wave systems a complete theory exists how to handle
the validity question of
NLS approximations and of generalizations of NLS approximations.
For proving the validity of NLS approximations for systems with quadratic nonlinearities
it is essential to consider the Fourier transformed systems.

{\bf i)} In Fourier space we have
$$
\partial_{t} \widehat{U}(k,t)
+  i k  \cA \widehat{U}(k,t)
+ \cE \widehat{U}(k,t)
=  \widehat{\cT}_2(\widehat{U} ,\widehat{U} )(k,t)  +
 \widehat{\cT}_3(\widehat{U} ,\widehat{U} ,\widehat{U} )(k,t)  .
$$

{\bf ii)} We diagonalize the linear part in Fourier space with $ \widehat{U}(k,t) = \widehat{S}(k) \widehat{V}(k,t) $, with $ \widehat{S}(k)  \in \C^{N \times N} $,
for $ k \in \R $. We assume for a moment that such a diagonalization is
possible for all $ k \in \R $.
After the diagonalization we find
\begin{equation} \label{sys5}
\partial_t \widehat{V}(k,t) = \widehat{D}(k)
\widehat{V}(k,t) +
\widehat{N}(\widehat{V})(k,t),
\end{equation}
where
\begin{eqnarray*}
\widehat{D}(k) & = & ( i \omega_n(k) )_{n \in \N}  = - (\widehat{S}(k))^{-1} ( i k  \cA
+ \cE) \widehat{S}(k) , \\
\widehat{N}(\widehat{V})(k,t) & = &
 \widehat{G}_2
(\widehat{V} ,\widehat{V} )(k,t)
+ \widehat{G}_3
(\widehat{V} ,\widehat{V} , \widehat{V})(k,t),
\end{eqnarray*}
with
\begin{eqnarray*}
\lefteqn{( \widehat{G}_2
(\widehat{V} ,\widehat{V} ))_j (k,t) }\\ & = &
 \sum_{j_1, j_2 =1}^N \int  \widehat{g}^j_{ j_1 j_2}(k,
k-k_1,k_1)
\widehat{V}_{j_1}(k-k_1,t) \widehat{V}_{j_2}(k_1,t) dk_1
\end{eqnarray*}
and
\begin{eqnarray*}
\lefteqn{( \widehat{G}_3
(\widehat{V} ,\widehat{V} ,\widehat{V} ))_j (k,t) } \\ & = &
 \sum_{j_1, j_2 ,j_3=1}^N \int \int \widehat{g}^j_{ j_1 j_2 j_3}(k,
k-k_1,k_1-k_2, k_2)
\\ && \qquad \qquad \times
\widehat{V}_{j_1}(k-k_1,t) \widehat{V}_{j_2}(k_1-k_2,t) \widehat{V}_{j_3}(k_2,t) dk_2 dk_1.
\end{eqnarray*}
Without loss of generality we assume that  $ \widehat{g}^j_{ j_1 j_2} $ is symmetric
in $ j_1 $ and $ j_2 $, as well as  $ \widehat{g}^j_{ j_1 j_2 j_3} $ is symmetric
in $ j_1 $, $ j_2 $ and $ j_3 $.
These are exactly the dispersive wave systems for which the
analysis  for the validity  of the NLS approximation has been carried out.
In the following we transfer the higher order validity results for systems with quadratic nonlinearities
and the validity results for solutions with more than one NLS-scaled wave-packet
to \eqref{sys1}.

For \eqref{sys1} we have
\begin{eqnarray*}
 \widehat{G}_2
(\widehat{V} ,\widehat{V} )(k,t) & = &
(\widehat{S}(k))^{-1} \widehat{\cT}_2 (\widehat{S} \widehat{V} ,\widehat{S}\widehat{V})(k,t) , \\
 \widehat{G}_3
(\widehat{V} ,\widehat{V} , \widehat{V})(k,t) & = &
(\widehat{S}(k))^{-1} \widehat{\cT}_3 (\widehat{S} \widehat{V} ,\widehat{S} \widehat{V} ,\widehat{S}\widehat{V})(k,t)
\end{eqnarray*}
If
$$
(\widehat{\cT}_2 ( \widehat{U} ,\widehat{U}))_j(k,t) =  \sum_{j_1, j_2 =1}^N  \int \widehat{b}^j_{ j_1 j_2}
\widehat{U}_{j_1}(k-k_1,t) \widehat{U}_{j_2}(k_1,t) dk_1
$$
with $ \widehat{b}^j_{ j_1 j_2} \in \R $ we have
$$
\widehat{g}^j_{ j_1 j_2}(k,
k-k_1,k_1) =   \sum_{j_3, j_4, j_5 =1}^N  ((\widehat{S}(k) )^{-1})_{j,j_3}  \widehat{b}^{j_3}_{ j_4 j_5} \widehat{S}_{j_4,j_1}(k-k_1) \widehat{S}_{j_5,j_2}(k_1)
$$
and similar for $ \widehat{\cT}_3 $ and $ \widehat{g}^j_{ j_1 j_2 j_3}(k,
k-k_1,k_1-k_2, k_2) $.
\medskip

{\bf iii)} Since the  derivation of the approximation equations is notationally more simple in physical space we transfer the diagonlized system \eqref{sys5}
back to physical space where for the components we have
\begin{equation} \label{sys6}
\partial_t {V}_n(x,t) = i  {\omega}_n(- i \partial_x)
{V}_n(x,t) +
{N}_{n}(V)(x,t)
\end{equation}
for $ n = 1,\ldots,N $ which is the starting point of the subsequent analysis.

\section{Higher order NLS approximation}
\label{sec3}

In this section first we recall the results from the existing literature
about the validity of higher order NLS approximations for \eqref{sys5}, see for instance \cite{Kal87,SU17}.
Secondly, we use these results to extend the results from
\cite{BJL24} with a pure cubic nonlinearity to the situation
with quadratic nonlinearities, i.e., $ \cT_2 \neq 0 $.
We consider so called polarized initial conditions which are of order
$ \mathcal{O}(\varepsilon) $ only in one component,
say, we make the NLS ansatz for $ V_{n_0} $ for an $ n_0 \in \{ 1,\ldots,N \} $.
The ansatz for the derivation of a higher order NLS approximation is then given by
\begin{eqnarray} \label{approx1}
V_{n_0}(x,t) = \sum_{m= - m_*}^{m_*} \sum_{j=0}^{j_*(m)}\varepsilon^{\beta_{n_0}(m)+ j} A_{n_0,m,j}(X,T) e^{im (k_0 x - \omega_0 t)}
\end{eqnarray}
and by
\begin{eqnarray} \label{approx2}
V_n(x,t) = \sum_{m= - m_*}^{m_*} \sum_{j=0}^{j_*(m)}\varepsilon^{\beta_n(m)+ j} A_{n,m,j}(X,T) e^{im (k_0 x - \omega_0 t)}
\end{eqnarray}
for $ n \in \{ 1,\ldots,N \}$ with $ n \neq n_0 $,
where $ X = \varepsilon (x-ct) $ and  $ T = \varepsilon^2 t $, and where $ m^* $ is a fixed chosen number and where $ j_*(m)  $ is defined below.
Herein, we have the basic temporal wave number $ \omega_0 = - \omega_{n_0}(k_0) $  and $ c = \frac{d}{dk}\omega_{n_0}(k_0) $
the linear group velocity.
The appearing numbers are given by
$$
\beta_{n_0}(m) = 1+ ||m|-1| , \qquad \beta_n(m) = 1+ ||m|-1| + 2 \delta_{|m|1},
$$
and
$$
j_*(m) = m_* - |m| - 2 \delta_{|m|1},
$$
with $ \delta_{mn} $ the Kronecker delta.
Plugging this ansatz into \eqref{sys6} and equating equal powers of $ \varepsilon $ and of
$  e^{im (k_0 x - \omega_0 t)} $ to zero gives
that $ A_{n_0,1,0} $ has  to satisfy an  NLS equation
\begin{equation} \label{NLS}
\partial_T A_{n_0,1,0}   =  i \nu_1 \partial_X^2 A_{n_0,1,0}  + i \nu_2 A_{n_0,1,0}  |A_{n_0,1,0} |^2,
\end{equation}
with coefficients $ \nu_1 $, $  \nu_2 \in \R $,
where $  A_{n_0,-1,0}  = \overline{A_{n_0,1,0}} $. We find
that the $ A_{n_0,1,j} $ and $ A_{n_0,-1,j} $ for $ j \geq 1 $ have to satisfy  linear inhomogeneous
Schrödinger equations and that all other $ A_{n,m,j} $  have to satisfy algebraic equations which can be solved w.r.t.
$ A_{n,m,j} $
if the non-resonance conditions
\begin{equation} \label{nonres1}
\omega_n(m k_0) - m  \omega_{n_0}(k_0) \neq 0
\end{equation}
are satisfied for all $ n \in \{ 1,\ldots,N \} $ and $ m \in \{ -m_*,\ldots,m_* \} $ except if $ (m,n) = (1,n_0) $.  The approximation constructed in this way
is called in the following $  \varepsilon \Psi_{n_0} $.

In case of no quadratic terms, i.e. $ \cT_2 = 0 $, this approximation can be justified
with an approximation theorem by a simple  application of Gronwall's inequality, cf. \cite{KSM92}.
In case of $ \cT_2 \neq 0 $ the quadratic terms have to be eliminated by some normal form
transformation.
This requires the validity of additional non-resonance conditions, namely
\begin{equation} \label{nonres2}
\inf_{n_1,n_2 \in \{ 1,\ldots,N \} }\inf_{k \in \R} |\omega_{n_1}(k) - \omega_{n_0}(k_0) - \omega_{n_2}(k-k_0) | > 0 .
\end{equation}
This condition can be weakened to
\begin{equation} \label{nonres3}
\sup_{n_1,n_2 \in \{ 1,\ldots,N \} }\sup_{k \in \R} |\frac{\widehat{g}^{n_1}_{ n_0 n_2}(k,
k_0,k-k_0)}{\omega_{n_1}(k) - \omega_{n_0}(k_0) - \omega_{n_2}(k-k_0)} | < \infty .
\end{equation}
Moreover, we assume that
\begin{equation} \label{assumeS}
\sup_{k \in \R} \| \widehat{S}(k) \| + \sup_{k \in \R} \| (\widehat{S}(k))^{-1} \|
< \infty.
\end{equation}
Then we have the following approximation theorem, cf. \cite{Kal87,Schn11OW,SU17}.
\begin{theorem} \label{theorem1}
For all $ m \in \N $ with $ m \geq 4 $  the
following holds. Assume the validity of \eqref{assumeS} and of the non-resonance conditions
\eqref{nonres1} and \eqref{nonres3}.
Let $ A_{n_0,1,0} \in C([0,T_0],H^{3 (m-3)+2}) $ be a solution of the NLS equation  \eqref{NLS}
and let $ \varepsilon \Psi_{n_0} $ be the approximation defined above with $ m^* = m - 1$.
Then there exist $ C > 0 $ and $ \varepsilon_0 > 0 $ such that for all
$ \varepsilon \in (0,  \varepsilon_0) $ we have solutions $ V $ of \eqref{sys6}
with
$$ \sup_{t \in [0,T_0/\varepsilon^2]} \sup_{x \in \R} | V(x,t) - \varepsilon \Psi_{n_0}(x,t) | \leq C \varepsilon^{m - 5/2}.
$$
\end{theorem}
\begin{remark}
{\rm
In case of $ \cT_2 = 0 $ from the non-resonance condition \eqref{nonres1}
the cases $ m \in \{-2,0,2\} $ only play a role for the higher order terms, i.e. $ j \geq 1 $, and \eqref{nonres3}
is no longer necessary.
}
\end{remark}
For the non-diagonalized system we obtain
\begin{corollary}
Under the assumptions of Theorem \ref{theorem1}
there exist $ C > 0 $ and $ \varepsilon_0 > 0 $ such that for all
$ \varepsilon \in (0,  \varepsilon_0) $ we have solutions $ U $ of \eqref{sys4}
with
$$ \sup_{t \in [0,T_0/\varepsilon^2]} \sup_{x \in \R} | U(x,t) - \varepsilon S \Psi_{n_0}(x,t) | \leq C \varepsilon^{m - 5/2}.
$$
\end{corollary}
For the original system \eqref{sys1} the higher order approximation result is as follows.
\begin{corollary} \label{cor1}
Under the assumptions of Theorem \ref{theorem1}
there exist $ C > 0 $ and $ \varepsilon_0 > 0 $ such that for all
$ \varepsilon \in (0,  \varepsilon_0) $ we have solutions $ \cU $ of \eqref{sys1}
with
$$ \sup_{\tau \in [0,T_0/\varepsilon]} \sup_{\xi \in \R} | \cU(\xi,\tau) -  S \Psi_{n_0}(\xi,\tau) | \leq C \varepsilon^{m - 7/2}.
$$
\end{corollary}
\begin{remark}
{\rm Theorem \ref{theorem1} should not be taken for granted since solutions
of order $ \mathcal{O}(\varepsilon) $ have to be controlled on an
$ \mathcal{O}(1/\varepsilon^2) $-time-scale.
}
\end{remark}
\begin{remark}
{\rm
The  previous approximation results, in particular Corollary \ref{cor1} guarantee that the NLS approximation
can be used for an effective simulation of
solutions of \eqref{sys1} to polarized initial conditions of the form
\eqref{init1}. As already said an initial condition is called polarized if in the diagonalized
system \eqref{sys6}
the initial condition is $ \mathcal{O}(\varepsilon) $ in only one component
at $ k = k_0 $.
The situation of initial conditions being $ \mathcal{O}(\varepsilon) $ in all components at $ k = k_0 $ is handled in the following Section \ref{sec4}.
}
\end{remark}
\begin{remark}
{\rm
For having an approximation of the form \eqref{approx1}-\eqref{approx2},
i.e., that $ A_{n_0,-1,0} $ also belongs to the $ V_{n_0} $-component
we need that $ \omega_{n_0}(k) = - \omega_{n_0}(-k) $ around
$ k =  k_0 $. In general this is only possible if $ \omega_{n_0} $ is defined with at least one jump. For the derivation of the amplitude equations
this jump should be chosen at a small but order $ \mathcal{O}(1) $ wave number   $ k \neq  0 $.
}
\end{remark}
\begin{remark}
{\rm The assumption on the diagonalization can be
weakened strong\-ly. Only a separation of the NLS modes
near $ k = \pm k_0 $ is necessary.
However, since this requires a complete rewriting of
all non-resonance conditions, cf. \cite{Kal87,Co02},
and gives less insight we prefer to stay at the chosen
presentation.
}
\end{remark}
\begin{remark}
{\rm
The eigenvalues for Example  \ref{example2} are the same as
for the
Klein-Gordon model \eqref{KG}. It is well known that
the non-resonance condition \eqref{nonres3} is satisfied for the
Klein-Gordon model \eqref{KG}. The matrix $ S $ can be defined in
Fourier space by
$$
\widehat{S}(k)  = \frac{1}{\sqrt{1+k^2}} \left( \begin{array}{cc}  1+ik & 1+ik  \\ -i \sqrt{1+k^2} &i \sqrt{1+k^2}   \end{array} \right) ,
$$
for which the validity of the assumption \eqref{assumeS}
is obvious since the limits exist for $ |k| \to \infty $.
%
%
%
%
%
%
}
\end{remark}

\section{The handling of more than one wave-packet}
\label{sec4}

It is the purpose of this section to explain how to handle more than one wave-packet, i.e,
how to handle the case of non-polarized initial conditions.
These are of order
$ \mathcal{O}(\varepsilon) $  in all components.
In order to handle this situation
we make an NLS approximation not only for $ V_{n_0} $ alone, but for all
$ V_n $ with $ n \in \{ 1,\ldots,N \} $. We restrict ourselves again to the case $ x \in \R $. Due to the fact that in general  the group velocities
$ \frac{d}{dk} \omega_n$
of the wave-packets are different, no consistent   ansatz of the  above form \eqref{approx1}-\eqref{approx2}
is possible.
However, by a slight modification of this ansatz a consistent ansatz
up to an error of order $ \mathcal{O}(\varepsilon^3)$ is possible.
In detail, the analysis made in \cite{CS12} for the interaction of two
NLS-scaled wave-packets for dispersive wave systems on the one hand can be
specialized and  on the other hand can be generalized
to handle the present  situation.

The ansatz of \cite{CS12} specialized and generalized to the present situation is given by
\begin{eqnarray*}
V_n(x,t) &= & \sum_{r=0}^2 \veps^{1+r} A_{n,1,r}(X_n,T)e^{iY_n} \\ && +
 \sum_{r=0}^2 \veps^{1+r} A_{n,-1,r}(X_n,T)e^{-iY_n} +
 M_{mixed,n},\\[.3cm]
X_n&=&X+\veps\omega_{n}'(k_0)t + \veps^2 \sum_{j \neq n}
\Psi_{n,j}^{(1)}(\X+\veps
\omega_{j}'(k_0)t,\T),\\[.3cm]
Y_n&=& k_0 x - \omega_{n}(k_0)t+ \sum_{l=1,2}\veps^l  \sum_{j \neq n}\Omega_{n,j}^{(l)}(X+\veps
\omega_{j}'(k_0)t,T),\\
X&=&\veps x, \quad T=\veps^2 t,
\end{eqnarray*}
with $  M_{mixed,n} = \mathcal{O}(\veps^2)$ determined below. It is a specialization in the
sense that all spatial wave numbers of the wave-packets are the same.
It is generalization in the
sense that more than two wave-packets are considered.

In case  of $ \cT_2 = 0 $ in \eqref{sys1} we recall the explicit formulas which determine
the amplitudes  $ A_{n,1,r} $, the phase shifts $ \Omega_{n,j}^{(1)} $,
the envelope shifts $ \Psi_{n,j}^{(1)} $, and
the second order corrections of the phase shifts and amplitudes
$ \Omega_{n,j}^{(2)} $.

Similar to \cite{CS12}, at $ \mathcal{O}(\varepsilon^3) $ we find that  $ A_{n,1,0} $
has to satisfy  the NLS equation
\begin{eqnarray} \label{NLS5}
\partial_2 A_{n,1,0}(X_n,T) & = & - i (\omega''_n(k_0)/2)
\partial_1^2 A_{n,1,0} (X_n,T)  \\ && \qquad+ 3
g^n_{nnn}(k_0,k_0,k_0,-k_0)|A_{n,1,0}(X_n,T) |^2
A_{n,1,0}(X_n,T) , \nonumber
\end{eqnarray}
and that  $ \Omega_{n,j}^{(1)} $ has to satisfy  the phase shift formula
$$
\Omega_{n,j}^{(1)}({X}_j,T)=\frac{6 g^n_{n j j}(k_0,k_0,k_0,-k_0)}{i(
\omega'_n(k_0) -\omega'_j(k_0) )}\int^{{X}_j}_{-\infty} |A_{j,1,0}(\zeta,T)|^2
d \zeta.
$$ In particular, both the NLS equation and the phase shift formula are analogous to (3.10) and (3.11) in \cite{CS12} respectively.
We have that  $\Omega_{n,j}^{(1)}$ is a real quantity because of $ g^n_{n j j}(k_0,k_0,k_0,-k_0) \in i \R $.  Therefore, it is a pure phase correction.
Moreover, we have a number of mixed terms which can be eliminated
by setting $ M_{mixed,n} $ equal to
$$
\sum_{r_1,r_2,r_3 = \pm 1}\sum_{j_1,j_2,j_3, r_1 j_1+ r_2 j_2 + r_3 j_3 \neq 1 }^N \veps^{3} M_{n,j_1,j_2,j_3}^{r_1,r_2,r_3}(X,T)
e^{i(r_1 Y_{j_1}+ r_2 Y_{j_2}+r_3 Y_{j_3} )} + c.c.,
$$
where
\begin{eqnarray*}
M_{n,j_1,j_2,j_3}^{r_1,r_2,r_3}(X,T) & = &
((\omega_{j_1}(r_1  k_0)+\omega_{j_2}(r_2  k_0)+  \omega_{j_3}(r_3  k_0)-
\omega_n((r_1+r_2+r_3) k_0))^{-1}  \\ && \qquad  \times
g^n_{j_1,j_2,j_3}((r_1+r_2+r_3)k_0,r_1 k_0,r_2 k_0,r_3 k_0)
\\ && \qquad \qquad \times
A_{j_1,r_1,0}(X_{j_1},T)A_{j_2,r_2,0}(X_{j_2},T)A_{j_3,r_3,0}(X_{j_3}
,T).
\end{eqnarray*}

 Comparable to (3.12) in \cite{CS12}, at $ \mathcal{O}(\varepsilon^4) $ we find that the $ A_{n,1,1} $
 solve  linear inhomogeneous
evolution equations
$$
\partial_2 A_{n,1,1} ({X}_1, T)= i (\omega''_1(k_0)/2)
\partial_1^2 A_{n,1,1}  ({X}_1, T)
+t_{n,1,1},
$$
where $ t_{n,1,1} $ is a function of $ A_{n,1,0} $ and $ A_{n,-1,0} $. Thus, $ A_{n,1,1} $ describes internal dynamics of a single pulse.

At $ \mathcal{O}(\varepsilon^4) $ we also find the envelope shift formula
$$
\Psi_{n,j}^{(1)}(\underline{X}_j,T)= C_{n,j}\int^{\underline{X}_j}_{-\infty}
|A_{j,1,0}(\zeta,T)|^2 d \zeta,
$$
with an  explicitly computable prefactor $C_{n,j}$, equivalent to (3.13) in \cite{CS12}.

The quantities $\Omega_{n,j}^{(2)}$ are determined at $ \mathcal{O}(\varepsilon^4) $.
The
real part is a second order correction to the phase shift, whereas its imaginary
part gives a correction to the amplitude. We refrain from explicitly displaying
the rather lengthy expression for $\Omega_{n,j}^{(2)}$ and only note that
it is pure integration of spatially localized terms similar to the expressions for
determining
$\Omega_{n,j}^{(1)}$ and $\Psi_{n,j}^{(1)}$.

Finally, there are even more mixed terms which we also do not display explicitly.

For the computation of the mixed terms again a number of
non-resonance conditions have to be satisfied.
For the computation of the quadratic mixed terms we need
\begin{equation} \label{nonres1a}
\omega_n((r_1+r_2) k_0) -   \omega_{j_1}(r_1 k_0) -
 \omega_{j_2}(r_2 k_0)  \neq 0
\end{equation}
for all $ n,j_1,j_2 \in \{ 1,\ldots,N \} $ and  $ r_1,r_2 \in \{-1,1\} $.
For the computation of the cubic mixed  terms we need
\begin{equation} \label{nonres1b}
\omega_n((r_1+r_2+r_3) k_0) -   \omega_{j_1}(r_1 k_0) -
 \omega_{j_2}(r_2 k_0) -   \omega_{j_3}(r_3 k_0) \neq 0
\end{equation}
for all $ n,j_1,j_2,j_3 \in \{ 1,\ldots,N \} $ and  $ r_1,r_2,r_3 \in \{-1,1\} $
with $ r_1 + r_2+ r_3 \not \in \{-1,1\} $.
For the computation of the quartic mixed  terms we need
\begin{equation} \label{nonres1c}
\omega_n((r_1+r_2+r_3+ r_4) k_0) -   \omega_{j_1}(r_1 k_0) -
 \omega_{j_2}(r_2 k_0) -   \omega_{j_3}(r_3 k_0) -    \omega_{j_4}(r_4 k_0)\neq 0
\end{equation}
for all $ n,j_1,j_2,j_3,j_4 \in \{ 1,\ldots,N \} $ and  $ r_1,r_2,r_3,r_4 \in \{-1,1\} $. For the computation of the phase shifts $ \Omega_{n,j}^{(r)} $ and the envelope shifts $ \Psi_{n,j}^{(1)} $
the additional condition
\begin{equation} \label{nonres1d}
\omega'_{n}(k_0) -
 \omega'_{j}(k_0) \neq 0,
\end{equation}
for all $ n,j \in \{ 1,\ldots,N \} $, with $ n \neq j $,
on the group velocities  is necessary.

In case of no quadratic terms, i.e. $ \cT_2 = 0 $, again  this approximation can be justified
with a simple  application of Gronwall's inequality, cf. \cite{KSM92}.
In case of $ \cT_2 \neq 0 $ the quadratic terms have to be eliminated by some normal form
transformation.
This requires  additional non-resonance conditions, cf. \cite[\S 11.5]{SU17}, namely
\begin{equation} \label{nonres2a}
\inf_{n_1,n_2,n_3 \in \{ 1,\ldots,N \} }\inf_{k \in \R} |\omega_{n_1}(k) - \omega_{n_2}(k_0) - \omega_{n_3}(k-k_0) | > 0 .
\end{equation}
This condition again can be weakened, namely  to
\begin{equation} \label{nonres3a}
\sup_{n_1,n_2,n_3 \in \{ 1,\ldots,N \} }\sup_{k \in \R} |\frac{\widehat{g}^{n_1}_{n_2 n_3}(k,
k_0,k-k_0)}{\omega_{n_1}(k) - \omega_{n_2}(k_0) - \omega_{n_3}(k-k_0)} | < \infty .
\end{equation}
\begin{remark}
{\rm
For the computation of the phase shifts $ \Omega_{n,j}^{(r)} $ and the envelope shift $ \Psi_{n,j}^{(1)} $
some integral has to be computed and so a certain spatial
localization  of the solutions of the NLS equations is necessary.
Therefore, we define the space $ H^s_m $
for $ s,m\in \N_0 $
as a subspace of $ H^s $
for which the norm
$$
\| A \|_{H^s_m} = \| A \rho \|_{H^s}
$$
is finite, where $ \rho(x) = (1+ x^2)^{m/2} $.
Local existence and uniqueness for the solutions of the NLS equations \eqref{NLS5} in spaces $ H^{s+m}_0 \cap H^s_m $  is well known if $ s \geq 1 $ and follows by application of  the variation of constant formula and using the fact that $ i \partial_X^2 $ is the generator of a strongly continuous
semigroup in $ H^{s+m}_0 \cap H^s_m $, cf. \cite{CKS95}.
}\end{remark}

By the above ansatz the residual terms formally are of order
$ \mathcal{O}(\varepsilon^5) $.
Therefore, we have the following approximation theorem, cf. \cite[Theorem 11.2.6]{SU17}.
\begin{theorem} \label{theorem2}
Assume the validity of \eqref{assumeS} and of the non-resonance conditions
\eqref{nonres1a}, \eqref{nonres1b}, \eqref{nonres1c},  \eqref{nonres1d},
and \eqref{nonres3a}.
Let the $ A_{n,1,0} \in C([0,T_0],H^{12}_2 \cap H^{14}_0) $ be
 solutions of the NLS equations  \eqref{NLS5}
and let $ \varepsilon \Psi $ be the approximation to these solutions defined above
in Section \ref{sec4}.
Then there exist $ C > 0 $ and $ \varepsilon_0 > 0 $ such that for all
$ \varepsilon \in (0,  \varepsilon_0) $ we have solutions $ V $ of \eqref{sys6}
with
$$ \sup_{t \in [0,T_0/\varepsilon^2]} \sup_{x \in \R} | V(x,t) - \varepsilon \Psi(x,t) | \leq C \varepsilon^{5/2}.
$$
\end{theorem}
For the non-diagonalized system we obtain
\begin{corollary}
Under the assumptions of Theorem \ref{theorem2}
there exist $ C > 0 $ and $ \varepsilon_0 > 0 $ such that for all
$ \varepsilon \in (0,  \varepsilon_0) $ we have solutions $ U $ of \eqref{sys4}
with
$$ \sup_{t \in [0,T_0/\varepsilon^2]} \sup_{x \in \R} | U(x,t) - \varepsilon S \Psi(x,t) | \leq C \varepsilon^{5/2}.
$$
\end{corollary}
For the original system \eqref{sys1} the  approximation result is as follows.
\begin{corollary} \label{cor2}
Under the assumptions of Theorem \ref{theorem2}
there exist $ C > 0 $ and $ \varepsilon_0 > 0 $ such that for all
$ \varepsilon \in (0,  \varepsilon_0) $ we have solutions $ \cU $ of \eqref{sys1}
with
$$ \sup_{\tau \in [0,T_0/\varepsilon]} \sup_{\xi \in \R} | \cU(\xi,\tau) -  S \Psi(\xi,\tau) | \leq C \varepsilon^{3/2}.
$$
\end{corollary}
\begin{remark}
{\rm As before, Theorem \ref{theorem2} should not be taken for granted since solutions
of order $ \mathcal{O}(\varepsilon) $ have to be controlled on an
$ \mathcal{O}(1/\varepsilon^2) $-time-scale.
}
\end{remark}
\begin{remark}
{\rm
The  previous approximation results, in particular Corollary \ref{cor2} guarantee that the NLS approximation
can be used for an effective simulation of
solutions of \eqref{sys1} to general  initial conditions of the form
\eqref{init1}. The polarization is no longer needed.
}
\end{remark}

\section{The higher dimensional situation}
\label{sec6}

In this section we explain how the results change if
$ x \in \R^d $ with $ d \geq 2 $ is considered.

The results from Section \ref{sec3} transfer in a straightforward way
from $ x \in \R $ to $ x \in \R^d $ with $ d \geq 2 $.
The group velocity in \eqref{approx1} and \eqref{approx2} is then given by $ c = \nabla_{k_0} \omega_0 \in \R^d $
where $ k \in \R^d $.
The NLS equation is then given by
\begin{equation} \label{NLSd}
\partial_T A_{n_0,1,0}   =  i \sum_{j_1, j_2}^d \nu_{j_1 j_2} \partial_{X_{j_1}}
\partial_{X_{j_2}} A_{n_0,1,0}
 + i \nu_2 A_{n_0,1,0}  |A_{n_0,1,0} |^2,
\end{equation}
where $ \nu_{j_1 j_2} = \frac12 \partial_{k_{j_1}}
\partial_{k_{j_2}} \omega_0|_{k= k_0} $.
The proof of the  approximation result is line for line the same.
In a similar way the results from Section \ref{sec4} can be modified.

However, in general the non-resonance condition to eliminate
the quadratic terms \eqref{nonres2} will not be valid in higher
space dimensions. In \cite{DHSZ16} it has been pointed out that
this problem can be solved by working in modulational
Gevrey spaces.

We would like to close the paper with the remark that
the results from Section \ref{sec4} are less relevant in higher space
dimensions due to the fact that for $ x \in \R^d $ the wave-packets have more space and so in general spatially localized solutions will miss each other.
In case that  the wave-packets are
not at the same place at the same time, there is a
decoupling up to high order.
The ansatz  is then given by
$$
V(x,t) =  \sum_{n_0=1}^N \varepsilon \Psi_{n_0}(x,t)
$$
where $ \varepsilon \Psi_{n_0} $ is the approximation defined
in \eqref{approx1}-\eqref{approx2} but now with
$$ X = X_{n_0}= \varepsilon (x - \nabla \omega_{n_0} |_{k = k_0} t)  + \varepsilon^{-1} X_{n_0,0}  = \widetilde{X} + \varepsilon^{-1}(X_{n_0,0} - \nabla \omega_{n_0} |_{k = k_0} T), $$
with $ \widetilde{X} = \varepsilon x $.
Suppose now that
\medskip

{\bf (ASS)} the lines $ (X,T) = ( X_{n_0,0} - \nabla \omega_{n_0} |_{k = k_0} T,T) \in \R^{d+1}$ have no intersection points,
i.e., assume that
the wave-packets are not at the same  place at the same time.

Then we have
\begin{theorem} \label{theorem1g}
For all $ m \in \N $ there exists a $ s \in \N $ such that the
following holds. Assume  the validity of \eqref{assumeS}, of the non-resonance conditions
\eqref{nonres1} and \eqref{nonres3}, and of the non-interaction assumption {\bf (ASS)}.
For all $ n_0 \in \{1,\ldots,N\} $
let $ A_{n_0,1,0} \in C([0,T_0],H^{s+ 3 m_*}_{m} \cap H^{s+ 3 m_*+m}_{0}) $ be solution of the NLS equations  \eqref{NLS}
and let the $ \varepsilon \Psi_{n_0} $ be the approximations defined
in \eqref{approx1}-\eqref{approx2} with the above modification.
Then there exist $ C > 0 $ and $ \varepsilon_0 > 0 $ such that for all
$ \varepsilon \in (0,  \varepsilon_0) $ we have solutions $ V $ of \eqref{sys6}
with
$$ \sup_{t \in [0,T_0/\varepsilon^2]} \sup_{x \in \R} | V(x,t) - \sum_{n_0=1}^N \varepsilon \Psi_{n_0}(x,t) | \leq C \varepsilon^{m - 1/2}.
$$
\end{theorem}
By the localization $ A_{n_0,1,0} \in C([0,T_0],H^{s+ 3 m_*}_{m} ) $
the interaction terms are of sufficiently high order w.r.t. $ \varepsilon $, cf. \cite{PW95,CCSU08}.

%
%
%
%
%
%
%
%
%
%
%

\bibliographystyle{alpha}
\bibliography{jahnke06}

\newcommand{\etalchar}[1]{$^{#1}$}
\begin{thebibliography}{CBCSU08}

\bibitem[BJL24]{BJL24}
Julian Baumstark, Tobias Jahnke, and Christian Lubich.
\newblock Polarized high-frequency wave propagation beyond the nonlinear
  {S}chr\"odinger approximation.
\newblock {\em SIAM J. Math. Anal.}, 56(1):454--473, 2024.

\bibitem[CBCSU08]{CCSU08}
Martina Chirilus-Bruckner, Christopher Chong, Guido Schneider, and Hannes
  Uecker.
\newblock Separation of internal and interaction dynamics for {NLS}-described
  wave packets with different carrier waves.
\newblock {\em J. Math. Anal. Appl.}, 347(1):304--314, 2008.

\bibitem[CBS12]{CS12}
Martina Chirilus-Bruckner and Guido Schneider.
\newblock Detection of standing pulses in periodic media by pulse interaction.
\newblock {\em J. Differential Equations}, 253(7):2161--2190, 2012.

\bibitem[CBS15]{SC15}
Martina Chirilus-Bruckner and Guido Schneider.
\newblock Interaction of oscillatory packets of water waves.
\newblock {\em Discrete Contin. Dyn. Syst.}, pages 267--275, 2015.

\bibitem[CBSU07]{CSU07}
Martina Chirilus-Bruckner, Guido Schneider, and Hannes Uecker.
\newblock On the interaction of {NLS}-described modulating pulses with
  different carrier waves.
\newblock {\em Math. Methods Appl. Sci.}, 30(15):1965--1978, 2007.

\bibitem[CKS95]{CKS95}
Walter Craig, Thomas Kappeler, and Walter Strauss.
\newblock Microlocal dispersive smoothing for the {S}chr\"odinger equation.
\newblock {\em Comm. Pure Appl. Math.}, 48(8):769--860, 1995.

\bibitem[CL09]{CL09}
Mathieu Colin and David Lannes.
\newblock Short pulses approximations in dispersive media.
\newblock {\em SIAM J. Math. Anal.}, 41(2):708--732, 2009.

\bibitem[Col02]{Co02}
Thierry Colin.
\newblock Rigorous derivation of the nonlinear {Schr{\"o}dinger} equation and
  {Davey}-{Stewartson} systems from quadratic hyperbolic systems.
\newblock {\em Asymptotic Anal.}, 31(1):69--91, 2002.

\bibitem[CW17]{CW18}
Patrick Cummings and C.~Eugene Wayne.
\newblock Modified energy functionals and the {NLS} approximation.
\newblock {\em Discrete Contin. Dyn. Syst.}, 37(3):1295--1321, 2017.

\bibitem[DH18]{DH18}
Wolf-Patrick D{\"u}ll and Max He\ss.
\newblock Existence of long time solutions and validity of the nonlinear
  {S}chr\"odinger approximation for a quasilinear dispersive equation.
\newblock {\em J. Differential Equations}, 264(4):2598--2632, 2018.

\bibitem[DHSZ16]{DHSZ16}
Wolf-Patrick D\"ull, Alina Hermann, Guido Schneider, and Dominik Zimmermann.
\newblock Justification of the 2{D} {NLS} equation for a fourth order nonlinear
  wave equation---quadratic resonances do not matter much in case of analytic
  initial conditions.
\newblock {\em J. Math. Anal. Appl.}, 436(2):847--867, 2016.

\bibitem[DLP{\etalchar{+}}11]{Schn11OW}
Willy D{\"o}rfler, Armin Lechleiter, Michael Plum, Guido Schneider, and
  Christian Wieners.
\newblock {\em Photonic crystals. {Mathematical} analysis and numerical
  approximation.}, volume~42 of {\em Oberwolfach Semin.}
\newblock Berlin: Springer, 2011.

\bibitem[DSW16]{DSW16}
Wolf-Patrick D\"ull, Guido Schneider, and C.~Eugene Wayne.
\newblock Justification of the nonlinear {S}chr\"odinger equation for the
  evolution of gravity driven 2{D} surface water waves in a canal of finite
  depth.
\newblock {\em Arch. Ration. Mech. Anal.}, 220(2):543--602, 2016.

\bibitem[D{\"u}l17]{D17}
Wolf-Patrick D{\"u}ll.
\newblock Justification of the nonlinear {S}chr\"odinger approximation for a
  quasilinear {K}lein-{G}ordon equation.
\newblock {\em Comm. Math. Phys.}, 355(3):1189--1207, 2017.

\bibitem[D{\"u}l21]{D21}
Wolf-Patrick D{\"u}ll.
\newblock Validity of the nonlinear {S}chr\"odinger approximation for the
  two-dimensional water wave problem with and without surface tension in the
  arc length formulation.
\newblock {\em Arch. Ration. Mech. Anal.}, 239(2):831--914, 2021.

\bibitem[Kal87]{Kal87}
L.~A. Kalyakin.
\newblock Asymptotic decay of a one-dimensional wave packet in a nonlinear
  dispersive medium.
\newblock {\em Mat. Sb. (N.S.)}, 132(174)(4):470--495, 592, 1987.

\bibitem[KSM92]{KSM92}
Pius Kirrmann, Guido Schneider, and Alexander Mielke.
\newblock The validity of modulation equations for extended systems with cubic
  nonlinearities.
\newblock {\em Proc. Roy. Soc. Edinburgh Sect. A}, 122(1-2):85--91, 1992.

\bibitem[MN02]{MN02}
Nader Masmoudi and Kenji Nakanishi.
\newblock From nonlinear {K}lein-{G}ordon equation to a system of coupled
  nonlinear {S}chr\"odinger equations.
\newblock {\em Math. Ann.}, 324(2):359--389, 2002.

\bibitem[PW95]{PW95}
R.~D. Pierce and C.~E. Wayne.
\newblock On the validity of mean-field amplitude equations for
  counterpropagating wavetrains.
\newblock {\em Nonlinearity}, 8(5):769--779, 1995.

\bibitem[Rau12]{Rauch12}
Jeffrey Rauch.
\newblock {\em Hyperbolic partial differential equations and geometric optics},
  volume 133 of {\em Graduate Studies in Mathematics}.
\newblock American Mathematical Society, Providence, RI, 2012.

\bibitem[Sch98]{Schn98}
Guido Schneider.
\newblock Justification of modulation equations for hyperbolic systems via
  normal forms.
\newblock {\em NoDEA Nonlinear Differential Equations Appl.}, 5(1):69--82,
  1998.

\bibitem[Sch05]{Schn05}
Guido Schneider.
\newblock Justification and failure of the nonlinear {S}chr\"odinger equation
  in case of non-trivial quadratic resonances.
\newblock {\em J. Differential Equations}, 216(2):354--386, 2005.

\bibitem[SU17]{SU17}
Guido Schneider and Hannes Uecker.
\newblock {\em Nonlinear {PDE}s}, volume 182 of {\em Graduate Studies in
  Mathematics}.
\newblock American Mathematical Society, Providence, RI, 2017.
\newblock A dynamical systems approach.

\bibitem[TW12]{TW12}
Nathan Totz and Sijue Wu.
\newblock A rigorous justification of the modulation approximation to the 2{D}
  full water wave problem.
\newblock {\em Comm. Math. Phys.}, 310(3):817--883, 2012.

\bibitem[Zak68]{Za68}
V.E. Zakharov.
\newblock Stability of periodic waves of finite amplitude on the surface of a
  deep fluid.
\newblock {\em Sov. Phys. J. Appl. Mech. Tech. Phys}, 4:190--194, 1968.

\end{thebibliography}

\end{document}